\documentclass{amsart}

\usepackage{amssymb} 
\allowdisplaybreaks
\usepackage{mathtools} 

\usepackage{enumitem}
\usepackage{booktabs}
\usepackage{multirow}
\usepackage{float}
\usepackage{hyperref}
\usepackage{algorithm} 
\usepackage[noend]{algorithmic} 

\newtheorem{theorem}{Theorem}[section]

\theoremstyle{definition}
\newtheorem{definition}[theorem]{Definition}
\newtheorem{example}[theorem]{Example}

\theoremstyle{remark}

\numberwithin{equation}{section}

\DeclareMathOperator{\ord}{ord}

\newcommand{\F}{\mathbb{F}}

\newcommand{\J}{\mathcal{J}}

\title{On pairs of primes with small order reciprocity}

\author{Craig Costello}
\address{Queensland University of Technology, Brisbane, Australia}
\email{craig.costello@qut.edu.au}

\author{Gaurish Korpal}
\address{University of Arizona, Tucson, United States}
\email{gkorpal@arizona.edu}
\thanks{Part of this work was done while Gaurish was an intern at Microsoft Research.}

\date{}

\begin{document}

\begin{abstract}
We give a sieving algorithm for finding pairs of primes with small multiplicative orders modulo each other. This problem is a necessary condition for obtaining constructions of $2$-cycles of pairing-friendly curves, which have found use in cryptographic applications. 
Our database of examples suggests that, except for a well-known infinite family of such primes, instances become increasingly rare as the size of the primes increase. This leads to some interesting open questions for which we hope our database prompts further investigation. 
\end{abstract}
\maketitle

\section{The problem}

This paper is concerned with finding pairs of primes that satisfy certain instances of the following definition. 

\begin{definition}[Primes with $(k,k')$-order reciprocity]\label{def:main}
    The prime numbers $p$ and $q$ have \emph{$(k,k')$-order reciprocity} if $\ord_q(p)=k$ and  $\ord_p(q)=k'$. That is, $k$ and $k'$ are the smallest natural numbers such that $p^k \equiv 1 \pmod q$ and $q^{k'}\equiv 1\pmod p$.
\end{definition}

Herein we are interested in instances of Definition~\ref{def:main} where $p$ and $q$ are \emph{large} and where $k$ and $k'$ are \emph{small}. To be more concrete about what we mean by large and small, we must first discuss why this problem finds practical relevance in cryptography. 
In 2014, Ben{-}Sasson, Chiesa, Tromer, and Virza~\cite{scalable} showed that certain types of pairing-based zero-knowledge proof systems can be instantiated in a very efficient and scalable way if we have a \emph{$2$-cycle} of \emph{pairing-friendly} elliptic curves. Let $p$ and $q$ be large primes. We say that the elliptic curves $E/\F_p$ and $E'/\F_q$ are a $2$-cycle if and only if $\#E(\F_p)=q$ and $\#E'(\F_q)=p$. Furthermore, we say that $E/\F_p$ is pairing-friendly if its \emph{embedding degree} $k=\ord_q(p)$ is small enough for the order-$q$ Weil and/or Tate pairings to be efficiently computable. Similarly, we say that $E'/\F_q$ is pairing-friendly $p$ if its embedding degree $k'=\ord_p(q)$ is suitably small. These embedding degrees need to be small enough that the pairings in $\F_{p^k}^\times$ and $\F_{q^{k'}}^\times$ can be computed efficiently. In practice, and indeed in this paper, we follow the standard rule-of-thumb and impose that both $k$ and $k'$ are no larger than $50$~\cite{taxonomy}. 

On the other hand, for the associated proof system to be cryptographically secure, the elliptic curve discrete logarithm problems (ECDLPs) and the finite field discrete logarithm problems (DLPs) all need to be computationally infeasible. In practice, this means $p$ and $q$ must both be at least $2^{160}$ for ECDLP security, and $k$ and $k'$ must be large enough for the fields $\F_{p^k}$ and $\F_{q^{k'}}$ to offer comparable DLP security. 

\subsection{A necessary condition for 2-cycles}\label{sec:2cycle}

The only known instances of $2$-cycles of ordinary pairing-friendly elliptic curves arise from the Miyaji-Nakabayashi-Takano (MNT) family~\cite{MNT}, which is parameterised as
\begin{align}\label{eq:MNT}
p = x^2-x+1 \quad\quad {\rm and} \quad\quad q=x^2+1,
\end{align}
from which it follows that $(k,k')=(\ord_q(p),\ord_p(q))= (4,6)$.  

It turns out that the MNT $2$-cycle is suboptimal in practical applications, because the embedding degrees are too small to optimally balance the ECDLP and DLP securities. With $(k,k')=(4,6)$, meeting the requisite DLP security requires that $\F_{p^4}$ and $\F_{q^6}$ are at least a few thousand bits~\cite[Table 1]{ScottG18}, which forces the sizes of $p$ and $q$ to be orders of magnitude larger than they would be in the optimal scenario where $k$ and $k'$ are larger. The subsequent hunt for $2$-cycles of ordinary elliptic curves with larger embedding degrees has only produced negative results~\cite{cycles-math,cycles-family}. For example, Chiesa, Chua and Weidner~\cite{cycles-math} proved that $2$-cycles of ordinary elliptic curves $E/\F_p$ and $E'/\F_q$ with $(k,k') \in  \{(5,10), (8,8), (12,12)\}$ do not exist. In what follows, we use a toy example to motivate the work in this paper by showing the existence of a $2$-cycle with $(k,k')=(12,12)$ under relaxed conditions, in light of recent work in~\cite{CCN24}.

\begin{example}\label{eg:1212}
The prime pair $(p,q) = (620461,15493)$ is such that $(k,k')=(\ord_q(p),\ord_p(q))=(12,12)$. Proposition 3 of~\cite{cycles-math} rules out the existence of a $2$-cycle with $q=\#E(\F_p)$ and $p=\#E'(\F_q)$. Indeed, the Hasse bound for $E'/\F_q$ in this instance is $15246 \leq \#E'(\F_q) \leq 15742$, ruling out $p=\#E'$. 
If we relax the $2$-cycle condition to allow cofactors as in~\cite[\S 7]{cycles-math} and to allow extension fields as in~\cite{CCN24}, we might instead seek to find two elliptic curves $E/\F_p$ and $E'/\F_{q^2}$ such that $\#E(\F_p) = hq$ and $\#E'(\F_{q^2})=h'p$. While such an $E$ does exist (with $h=40$), there is no multiple of $p$ inside the Hasse interval for $E'/\F_{q^2}$; $h'=386$ is too small for $h'p$ to be inside the interval, while $h'=387$ is too large. 
If, however, we instead consider a genus 2 curve $C/\F_q$, close examination of the possible group orders of its Jacobian (see~\cite[\S 2.1]{GaudryKS11}) reveals that $h' \in \{374, \dots, 399\}$ are all plausible cofactors. Indeed, we readily find one such example as $C/\F_q \colon y^2=x^6 +6611x^5 + 13858x^4 + 6818x^3 +5652x^2 + 10423x + 1795$, which is such that $\#\J_C(\F_q)=383 \cdot p$, forming a $2$-cycle with $E/\F_p \colon y^2=x^3+30984x+426966$, which is such that $\#E(\F_p)=40 \cdot q$. Both $\J_C$ and $E$ are ordinary and have embedding degree 12.
\end{example}

The above example shows that it is possible to construct ordinary $2$-cycles under relaxed conditions, even when their existence has been ruled out in the case where both curves are elliptic curves of prime order, so long as we have primes with $(k,k')$-order reciprocity\footnote{The notion and terminology \emph{order reciprocity} is found in a 2019 post by Lee~\cite{SLee}.}. In our view, the existence of large primes with small order reciprocity is \emph{the} fundamental and necessary requirement for determining whether relaxed $2$-cycles can exist in the interesting cryptographic ranges. 

\subsection{Open questions}

One can show that there are no prime pairs with $(2,2)$-order reciprocity, but infinitely many prime pairs with $(p-1,q-1)$-order reciprocity\footnote{For the former, $p<q$ implies $p=q-1$, the only option for which is $(p,q)=(2,3)$, which has $\ord_p(q)=1$. For the latter, see~\cite{FredH}.}. However, in between these two extremes there are a number of questions that are of potential relevance to cryptography. For example, we saw above that the pair $(p,q) = (620461,15493)$ has $(12,12)$-order reciprocity, but this is the only such example we found. We state some open questions that warrant further investigation: 
\begin{enumerate}
\item \emph{Is $(620461,15493)$ the only prime pair with $(12,12)$-order reciprocity?}
\item \emph{Are there any fixed values of $(k,k')$ with ${\rm min}(k,k') > 4$ for which there are an infinite number of primes with $(k,k')$-order reciprocity?} 
\item \emph{Are there any fixed values of $(k,k')$ with ${\rm min}(k,k') > 2$ for which there are no pairs of primes with $(k,k')$-order reciprocity?} 
\end{enumerate}
Finally, we pose the question of interest that underlies this work: 
\begin{enumerate}[resume]
\item {\it Are there {\bf any} large pairs of primes with small $(k,k')$-order reciprocity?}
\end{enumerate}

In the search for answers to the above questions, we derived an algorithm for finding primes with $(k,k')$-order reciprocity -- this is presented in Section~\ref{sec:algo}. We used this algorithm to exhaustively search up to the 200 millionth prime for pairs with small order reciprocity  -- the database of examples we found is summarised in Section~\ref{sec:data}. 

We reiterate that we are not interested in the known examples which come from the MNT parameterisation in~\eqref{eq:MNT}, and that in practice we are therefore interested in examples with, say, $5 \leq k, k' \leq 50$. By ${\it large}$ primes, we mean cryptographically large, i.e. at least $2^{160}$. 

\section{The algorithms}\label{sec:algo}

Let $p$ be a prime number and let $\alpha_p \in \mathbb{F}_p$ be such that $\mathbb{F}_p^\times = \langle \alpha_p\rangle$, where elements in $\mathbb{F}_p^\times$ are represented by the corresponding least positive remainder in $(\mathbb{Z}/p\mathbb{Z})^\times$. Then Algorithm~\ref{alg:reciprocity} returns all $q<p$ such that $\ord_q(p)=k$ and $\ord_p(q)=k'$. Furthermore, let $P$ be a list of prime numbers and $A$ be a list such that $A[i]$ is a primitive element modulo $P[i]$. Algorithm~\ref{alg:primepairs} finds all the primes pairs $(p,q)$ with $p\in P$ and $q<p$ such that $\ord_q(p) = k$ and $\ord_p(q) = k'$ or $\ord_q(p) = k'$ and $\ord_p(q) = k$.

\begin{algorithm}[!ht]
    \caption{$\mathtt{reciprocity}(p,\alpha_p, k,k')$}\label{alg:reciprocity}
    \begin{algorithmic}[1]
    \STATE $Q\leftarrow \texttt{[]}$ \COMMENT{empty list}
    \IF{$p \equiv 1\pmod {k'}$}
        \STATE $t \leftarrow \frac{p - 1}{k'}$
        \FOR{$s = 1$ \TO $k' - 1$}\label{line:for}
            \STATE $m \leftarrow s \cdot  t$\label{line:mult}
            \IF{$\gcd(m, p - 1) = t$}\label{line:GCD}
                \STATE $q \leftarrow \alpha_p^m$ \label{line:alphap}
                \IF{$q \equiv 1 \pmod k$ \AND $q$ is prime}\label{line:primality}
                    \STATE // we have found $q$ such that $\ord_p(q)=k'$.
                    \IF{$p^k\equiv 1 \pmod q$}
                        \STATE $b\leftarrow 1$
                        \FOR{$d = 1$ \TO $k - 1$}
                            \IF{$k \equiv 0 \pmod d$  \AND $p^d \equiv 1 \pmod q$}\label{line:if}
                                    \STATE $b\leftarrow 0$
                                    \STATE break \COMMENT{$\ord_q(p)\neq k$}
                                \ENDIF
                        \ENDFOR
                        \IF{$b = 1$}
                            \STATE $Q\text{.append}(q)$ \COMMENT{$\ord_q(p) = k$}
                        \ENDIF
                    \ENDIF
                \ENDIF
            \ENDIF
        \ENDFOR
    \ENDIF
    \RETURN $Q$
    \end{algorithmic}
\end{algorithm}

\begin{algorithm}[!ht]
    \caption{$\mathtt{primepairs}(P,A, k,k')$}\label{alg:primepairs}
    \begin{algorithmic}[1]
        \STATE $L\leftarrow \texttt{[]}$ \COMMENT{empty list}
        \FOR{$p$ in $P$}
            \STATE $Q\leftarrow \mathtt{reciprocity}(p,\alpha_p, k,k')$
            \IF{$\#Q \neq 0$}
                \FOR{$q\in Q$}
                    \STATE $L\text{.append}((p,q))$
                \ENDFOR
            \ENDIF
            \IF{$k\neq k'$}
                \STATE $Q'\leftarrow \mathtt{reciprocity}(p,\alpha_p, k',k)$
                \IF{$\#Q' \neq 0$}
                    \FOR{$q'\in Q'$}
                        \STATE $L\text{.append}((q',p))$
                    \ENDFOR
                \ENDIF
            \ENDIF
        \ENDFOR
        \RETURN $L$
    \end{algorithmic}
\end{algorithm}

\subsection{Proof of correctness} We start by showing that Algorithm~\ref{alg:reciprocity} returns all primes $q<p$ such that $\ord_q(p) = k$ and $\ord_p(q) = k'$. Any such $q$ must be the least positive remainder of some power of $\alpha_p^t$ with $t=(p-1)/k'$ in $(\mathbb{Z}/p\mathbb{Z})^\times$; Lines~\ref{line:for}-\ref{line:mult} exhaust these possibilities, while Line~\ref{line:GCD} avoids testing those $q$ which cannot have exact order $k'$. Line~\ref{line:if} only keeps those $q$ which are prime and for which elements of order $k$ exist in $(\mathbb{Z}/q\mathbb{Z})^\times$. The remainder of Algorithm~\ref{alg:reciprocity} simply keeps those $p$ which have exact order $k$ in $(\mathbb{Z}/q\mathbb{Z})^\times$. Algorithm~\ref{alg:primepairs} loops over all combinations of $p, q \in P$ and, for $k \neq k'$, makes an additional call to Algorithm~\ref{alg:reciprocity} to reverse the order of $k$ and $k'$.\\

\subsection{Complexity analysis} Let $p$ be the largest prime in the list $P$, whose cardinality is $N$. We will analyse the time complexity of Algorithm~\ref{alg:primepairs} as $p \rightarrow \infty$, assuming that $k$ and $k'$ remain fixed, small constants. We will also assume that, as $p \rightarrow \infty$, the average size of the primes in $P$ is in $O(p)$; this assumption holds if $P$ contains all primes within a given interval and, in particular, if $P$ contains all primes up to $p$. Algorithm~\ref{alg:primepairs} makes either one or two calls to Algorithm~\ref{alg:reciprocity} for each prime in $P$, so it follows that the time complexity of Algorithm~\ref{alg:primepairs} is in $O(N \cdot R)$, where $R$ is the complexity of Algorithm~\ref{alg:reciprocity} on input of a prime of size in $O(p)$. 

We will now proceed step by step through Algorithm~\ref{alg:reciprocity}, stating the time complexity of each step together with the probability of advancing through the control statements. These can be accumulated in a straightforward way to obtain a precise big-$O$ complexity of Algorithm~\ref{alg:reciprocity}, but this becomes a rather messy and nonillustrative expression, so we will eventually ignore all logarithmic factors to instead state a final big-$\tilde{O}$ time complexity. 

For small, fixed $k$ and $k'$, calls to Algorithm~\ref{alg:reciprocity} advance to Line~\ref{line:GCD} $O(p)$ times on average. The time complexity of the GCD algorithm here is in $O(\log{p})$ and, heuristically, $O(1)$ of these tests advance to Line~\ref{line:alphap}. The asymptotic complexity of this exponentiation is in $O(\log{p} \log^2{(\log{p})} \log{\log{\log{p}}})$, using the Sch{\"o}nhage--Strassen algorithm~\cite{SchoenhageStrassen1971}. In Line~\ref{line:primality}, an optimised probabilistic primality test (like Miller-Rabin~\cite{Rabin1980}) has time complexity $O(\log^3{p})$, while the deterministic AKS primality test~\cite{Agrawal2004} has time complexity $O(\log^6{p})$. Since the size of the prime $q$ is in $O(p)$, the probability of advancing to the remaining steps is in $O(1/\log{p})$ by the Prime Number Theorem.

It follows that Algorithm~\ref{alg:reciprocity} has time complexity $\tilde{O}(p)$, and thus that Algorithm~\ref{alg:primepairs} has complexity $\tilde{O}(N \cdot p)$. In particular, if $P$ contains all the primes up to $p$, then Algorithm~\ref{alg:primepairs} has complexity $\tilde{O}(p^2)$ by the Prime Number Theorem. 

\section{The database}\label{sec:data}
We used the \texttt{parallel GP} interface of \texttt{PARI/GP}~\cite{PARI2} running on AMD Ryzen Threadripper PRO 3995WX with 128 GB memory to perform all the computations. The \texttt{PARI/GP} implementation of the algorithms and the data are available in 

\begin{center}
    \url{https://github.com/gkorpal/order-reciprocity}.\\[1em]
\end{center}

Let $P_1$ be the list of first 100 million primes and $P_2$ be the list of the next 100 million primes. We first performed the primitive element pre-computation for $P_1$ and $P_2$ using the in-built \texttt{znprimroot} function \cite[Algorithm 1.4.4]{cohen}. 

We then ran Algorithm~\ref{alg:primepairs} to find prime pairs in $P_1$ and for prime pairs (where the largest prime is) in $P_2$ with $(k,k')$-order reciprocity, for $2\leq k \leq k' \leq 51$. We could not parallelise the outermost for loop of this algorithm because of the memory sharing restrictions. However, we utilized the \texttt{parallel GP} interface to simultaneously run the computations for $(k,k')\in \{\{m,m+1\}\times\{n,n+1\} \mid m,n\in 2\mathbb{Z} \ \text{and}\ 2\leq m\leq n\leq 50\}$ such that $k\leq k'$.

Tables~\ref{tab1} and~\ref{tab2} present sample data, while the repository offers access to the full dataset. Table~\ref{tab1} highlights that $(4,6)$-order reciprocity remains the most prevalent, even as prime size increases. Meanwhile, Table~\ref{tab2} shows that larger prime sizes result in most $(k, k')$-order reciprocities yielding no prime pairs.

\begin{table}[!ht]
    \caption{Counts of prime pairs with selected $(k,k')$-order reciprocity. The bottom row counts the number of pairs when the larger prime is in $P_2$.}
\label{tab1}
    \begin{tabular}{llcccccccc}
        \toprule
        \multicolumn{2}{l}{$(k,k')$} & $(4,6)$ & $(3,10)$ & $(3,14)$ & $(4,46)$ & $(2,35)$ & $(11,49)$ & $(15,45)$   & $(38,45)$ \\ \midrule
        \#prime & $P_1$ & 738 & 20 & 14 &  8 & 5 & 4 & 3 & 2 \\ 
 \ \  pairs & $P_2$ & 258 & 0 & 0 &  2 & 2 & 3 & 2 & 2 \\ 
        \bottomrule
    \end{tabular}
\end{table}

\begin{table}[!ht]
   \caption{Counts of $(k,k')$-order reciprocities with selected number of prime pairs. The bottom row counts the number of $(k,k')$-order reciprocities when the larger prime is in $P_2$.
}
\label{tab2}
   \begin{tabular}{llcccccccccr}
       \toprule
       \multicolumn{2}{c}{\#prime pairs}  & 0    & 1   & 2    & 3   & 4 to 12 & 14 & 20 & 258 & 738 & Total\\ \midrule
       \multirow{2}{*}{\# $(k,k')$}& $P_1$  & 97   & 147 & 209 & 231 & 588   & 1  & 1  & 0   & 1 & 1275 \\
       & $P_2$  & 1108 & 159 & 6   & 1   & 0     & 0  & 0  & 1   & 0 & 1275 \\ \bottomrule
   \end{tabular}
\end{table}

\bibliographystyle{amsalpha}
\bibliography{refs}

\end{document}